\documentclass[12pt]{amsart}
\usepackage{amsmath,amsfonts,euscript,amscd,amsthm,amssymb,upref,graphics}

\theoremstyle{plain}
\newtheorem{theorem}{Theorem}
\swapnumbers

\newtheorem{lemma}[subsection]{Lemma}

\theoremstyle{definition}
\newtheorem{definition}[subsection]{Definition}
\newtheorem{example}[subsection]{Example}

\newtheorem{nothing*}[subsection]{}

\newcommand{\rien}[1]{}

\newcommand{\Aut}{ \operatorname{{\rm Aut}}}

\newcommand{\C}{\ensuremath{\mathbb{C}}}
\newcommand{\Reals}{\ensuremath{\mathbb{R}}}

\newcommand{\emb}{\hookrightarrow}

\renewcommand{\epsilon}{\varepsilon}
\renewcommand{\phi}{\varphi}
\renewcommand{\emptyset}{\varnothing}
\begin{document}
\renewcommand{\baselinestretch}{1.07}

\title[Proper disks in $\mathbb C^2$]
{Proper holomorphic disks in the complement of varieties in $\mathbb
C^2$}

\author{Stefan Borell}
\address{Mathematisches Institut \\ Universit\"at Bern
    \\Sidlerstr. 5
   \\ CH-3012 Bern, Switzerland}
\email{borell@math.unibe.ch}

\author{Frank Kutzschebauch}
\address{Mathematisches Institut \\ Universit\"at Bern
    \\Sidlerstr. 5
   \\ CH-3012 Bern, Switzerland}
\email{Frank.Kutzschebauch@math.unibe.ch}

\author{Erlend Forn\ae ss Wold}
\address{Mathematisches Institut \\ Universit\"at Bern
    \\Sidlerstr. 5
   \\ CH-3012 Bern, Switzerland}
\email{erlend.wold@math.unibe.ch}

\thanks{All authors supported by Schweizerische Nationalfonds grant 200021-116165/1}


\begin{abstract}
For any analytic set $X\subset\mathbb C^2$ there exists a proper
holomorphic embedding $\phi:\triangle\emb\mathbb C^2$  such that
$\phi(\triangle)\cap X=\emptyset$.
\end{abstract}
\maketitle \vfuzz=2pt

\vfuzz=2pt

\section{Introduction}

In \cite{fg} the authors proved that there are proper holomorphic
disks in $\mathbb C^2$ that avoid the set $\{zw=0\}$. Furthermore
they said that it would be interesting to know whether there could
be such disks avoiding any finite set of complex lines. We show that
there exist properly \emph{embedded} disks satisfying the more
general property:

\begin{theorem}\label{main}
Let $X$ be any closed proper subvariety of $\mathbb C^2$ and let
$\triangle$ denote the unit disk in $\mathbb C$.  Then there exists
a proper holomorphic embedding $\phi:\triangle\emb\mathbb C^2$ such
that $\phi(\triangle)\cap X=\emptyset$.
\end{theorem}

It should be remarked that the corresponding result for $\C$ instead of the disk is false since
Kobayashi hyperbolicity of $\C^2\setminus X$ is an obstruction. In fact it is known that any analytic subset $X$ of $\C^2$ can be embedded in a different way $f: X \emb \C^2$ into $\C^2$ such that
the complement $\C^2 \setminus f(X)$ is Kobayashi hyperbolic (see \cite{BF}, \cite{BK}). In such
a situation there is not even a nonconstant holomorphic map from $\C$ into that complement.
An easier example is the following:

\begin{example} Let $X$ be the union of the following three lines in $\C^2_{z,w}$:
$$ l_1 = \{w=0\}, \quad l_2 = \{w=1\}, \quad l_3 = \{z=w\} $$
If a   holomorphic map $\varphi : \C \to \C^2$ avoids $l_1$ and $l_2$ it is of the form
$\varphi (\theta) = ( f(\theta), c)$ since the projection $\pi_w \circ \varphi $ is a map from $\C$ into
$\C\setminus \{0, 1\}$ and thus constant $=c$. For $\varphi$ to be an embedding means that $f(\theta ) =
a \theta + b, a \ne 0$. Therefore the image of $\varphi$ meets $l_3$. 

Note that in this example $\C^2 \setminus X$ is not Kobayashi hyperbolic. The maps $\theta \mapsto
(\exp{\theta} + c , c)$ provide nondegenerate holomorphic maps from $\C$ into $\C^2\setminus X$ if $c\neq 0,1$.
\end{example}

Hyperbolicity of $\C^2 \setminus X$ is the reason why additional interpolation on discrete (or even finite)
sets is not possible in general for embeddings as in our theorem.

\section{Construction}

Recall the following (simplified) definition from \cite{KLW}:

\begin{definition}
Given a smooth real curve in $\C^2$ without self intersection
$\Gamma = \{\gamma (t) : t \in [0, \infty )$ or $ t \in
(-\infty,\infty )\} $. We say that $\Gamma$ has the {\sl nice
projection property} if there is a holomorphic automorphism $\alpha
\in \Aut_{hol} (\C^2) $ of $\C^2$ such that, if $\beta (t)=\alpha
(\gamma (t))$, $\Gamma^\prime = \alpha (\Gamma)$ and $\pi_1:\C^2\to \C$ denotes the projection onto the first coordinate, then the following
holds:

\begin{itemize}

\item[(1)] $\lim_{|t|\to \infty }|\pi_1 (\beta (t))| =\infty$,
\item[(2)] There is an $M \in\Reals$ such that for all $R\ge M$ we have that $\C\setminus (\pi_1
(\Gamma^\prime)) \cup \overline\Delta_R)$ does not contain any
relatively compact connected components.
\end{itemize}

\end{definition}

\medskip

Let $W$ denote the set $W:=\overline\triangle\setminus\{1\}$  and
let $\Gamma$ denote the set $\Gamma:=\{z\in W;|z|=1\}$.  We will say
that a subset $\tilde W\subset W$ is \emph{$\mathfrak{b}$-nice} if $\tilde
W$ has a smooth boundary, and if there is a disk $D$ centered at
$\{1\}$ such that $W\cap D=\tilde W\cap D$. We let $\tilde\Gamma$
denote $\partial\tilde W\cap W$.  Note that if $\phi(W)$ is an
embedding such that $\phi(\Gamma)$ has the nice projection property,
then $\phi(\tilde\Gamma)$ has the nice projection property.  This is
because the two embedded curves are the same near infinity.

\medskip

The following lemma will provide us with the inductive step in our
construction:

\begin{lemma}\label{lemma}
Let $X$ be a closed subvariety of $\mathbb C^2$ and let
$\phi:W\emb\mathbb C^2$ be a smooth embedding, holomorphic on the
interior, such that

\medskip

$(i)$ $\mathrm{lim}_{j\rightarrow\infty}\|\phi(z_j)\|=\infty$ for
all $\{z_j\}\subset W$  with $z_j\rightarrow 1$, \

$(ii)$ $\phi(\Gamma)$ has the nice projection property, \

$(iii)$ $\phi(W)\cap X\cap\overline{\mathbb B}_N=\emptyset$ \
($N\in\mathbb N)$,

$(iv)$ $\phi(\Gamma)\subset\mathbb C^2\setminus\overline{\mathbb
B}_{N+1}$.

$(v)$ $\phi(W)$ intersects $\partial\mathbb B_N$ transversally.

\medskip

Let $S$ denote the set $S:=\phi^{-1}(\phi(W)\cap\overline{\mathbb
B}_N)$, let $V$ be a connected component of $S$ and let
$\epsilon>0$. Then there exists a $\mathfrak{b}$-nice subset $\tilde
W$ of $W$ with $V\subset\subset\tilde W$, $\tilde W\cap(S\setminus
V)=\emptyset$, and an embedding $\tilde\phi:\tilde W\emb\mathbb C^2$
(smooth, and holomorphic on the interior) such that

\medskip

$(a)$ $\mathrm{lim}_{j\rightarrow\infty}\|\tilde\phi(z_j)\|=\infty$
for all $\{z_j\}\subset\tilde W$  with $z_j\rightarrow 1$, \

$(b)$ $\tilde\phi(\tilde\Gamma)$ has the nice projection property, \

$(c)$ $\tilde\phi(\tilde W)\cap
X\cap\overline{\mathbb{B}}_{N+1}=\emptyset$, \

$(d)$ $\tilde\phi(\tilde\Gamma)\subset\mathbb
C^2\setminus\overline{\mathbb B}_{N+2}$, \

$(e)$ $\tilde\phi(\tilde W)$ intersects $\partial\mathbb B_{N+1}$
transversally, \

$(f)$ $\|\tilde\phi - \phi\|_V<\epsilon$, \

$(g)$ $\tilde\phi(\tilde W\setminus V)\subset\mathbb
C^2\setminus\mathbb{B}_{N-\epsilon}$.
\end{lemma}

\begin{proof}
It is not hard to see that $(i)$ and $(iv)$ implies that there exist
positive real numbers $0<r,\delta<1$  such that the set $A_r:=\{z\in
W;|z|\geq r\}$ satisfies
$$
(*) \ \phi(A_r)\subset\mathbb C^2\setminus\overline{\mathbb
B}_{N+1+\delta}.
$$
It follows that the set $P:=\phi^{-1}(\phi(W)\cap
(X\cap\overline{\mathbb B}_{N+1}))$ is a finite set of points since
the total intersection set $Z:=\phi^{-1}(\phi(W)\cap X)$ is analytic
in $W$ and since $P\subset Z\cap(W\setminus A_r)$.  The set $S$ is
clearly also contained in $W\setminus A_r$ and by $(v)$ it is a
finite disjoint union of smoothly bounded sets. \

Now for an arbitrarily small neighborhood $\mathcal{N}$ of $P$  we
have that $\mathrm{dist}(\phi((\overline{W\setminus
A_r})\setminus\mathcal N),X\cap\overline{\mathbb B}_{N+1})>0$.  Since
we also have $(*)$ we get that
$$
\mathrm{dist}(\phi(W\setminus\mathcal N),X\cap\overline{\mathbb
B}_{N+1})>0.
$$

This means that we may choose a $\mathfrak{b}$-nice domain $\tilde
W\subset W\setminus (P\cup(S\setminus V))$ such that
$V\subset\subset\tilde W$ and such that $\mathrm{dist}(\phi(\tilde
W),X\cap\overline{\mathbb B}_{N+1})>0$. Note that
$\phi(\tilde\Gamma)\subset\mathbb
C^2\setminus\overline{\mathbb{B}}_{N}$ and that $\phi(\tilde\Gamma)$
has the nice projection property since $\tilde\Gamma$ is the same as
$\Gamma$ near 1.  \

Since $K:=\overline{\mathbb B}_{N}\cup(X\cap\overline{\mathbb
B}_{N+1})$ is polynomially convex there is an open neighborhood
$\Omega$ of $K$ such that $\overline\Omega$ is polynomially convex
and such that $\overline\Omega\cap\phi(\tilde\Gamma)=\emptyset$.
Thus by Lemma 2.3 in \cite{KLW} there exists a $\Phi \in \Aut_{hol}
(\C^2)$ such that $\|\Phi-Id\|_{\overline\Omega}<\epsilon$ and such
that $\Phi(\phi(\tilde\Gamma))\subset\mathbb
C^2\setminus{\overline{\mathbb B}_{N+2}}$.  By possibly having to
decrease $\epsilon$ we may assume that $\Phi(\phi(\tilde W))\cap
(X\cap\overline{\mathbb B}_{N+1})=\emptyset$  and so we may put
$\tilde\phi:=\Phi\circ\phi$.  The conditions $(a),(c),(d)$ and $(f)$
are then immediate.  Since $\Phi$ is an automorphism we have that
$\tilde\phi(\tilde\Gamma)$ has the nice projection property and so
we get $(b)$.  Condition $(g)$ follows since we chose $\tilde W$
such that $\phi(\tilde W\setminus V)\subset\mathbb
C^2\setminus\mathbb B_N$ and because
$\|\Phi-Id\|_{\overline\Omega}<\epsilon$.  Finally, consider the
case where the intersection of $\tilde\phi(\tilde W)$  with any
sphere $\partial\mathbb B_\rho$ is not transversal.  In that case
there is a point $z\in\tilde W$ with $\|\tilde\phi(z)\|=\rho$ and
$\langle\tilde\phi(z),\mathrm{d}\tilde\phi(z)\rangle=0$.  Hence the
set of problematic points is analytic and thus discrete in $\tilde
W.$  So there exist $\rho$'s arbitrarily close to $1$  such that the
intersection of $\tilde\phi(\tilde W)$  with $\partial\mathbb
B_{\rho(N+1)}$ is transversal.  Thus there are arbitrarily small
linear perturbations of $\tilde\phi$ that give us $(e)$, and the
other properties are clearly preserved.

\end{proof}

\emph{Proof of Theorem \ref{main}}:  We will inductively construct
an increasing sequence of simply connected sets in the unit disk
along with a corresponding sequence of holomorphic embeddings.

\medskip

To start the induction we embed the disk into $\mathbb C^2$ as
follows: Start by letting $f_1:\overline\triangle\emb\mathbb C^2$ be
the map $z\mapsto (3z,0)$.  We may of course assume that
$f(\overline\triangle)\cap X=\emptyset$.  For $\delta>0$ let
$f_\delta$ denote the rational map $f_\delta:\mathbb
C^2\rightarrow\mathbb C^2$ given by $(z,w)\mapsto
(z,w+\frac{\delta}{z-1})$.  Let $W$ be as in Lemma \ref{lemma}.  If
we put $\phi_1:=f_\delta\circ f_1:W\rightarrow\mathbb C^2$ for a
small enough $\delta$ it is not hard to verify that all conditions
in Lemma \ref{lemma} are satisfied with $N=1$ (to get the nice
projection property, project to the $w$-axis).  Let
$U_1:=\phi_1(W)\cap\overline{\mathbb{B}}_1$ - a set we may assume to
be connected and (automatically) simply connected - and choose
$\epsilon_1>0$ such that if $\psi:\overline U_1\rightarrow\mathbb
C^2$ is any holomorphic map with $\|\psi-\phi_1\|_{\overline
U_1}<\epsilon_1$ then $\psi$ is an embedding and $\psi(\overline
U_1)\cap X=\emptyset$.  Choose $\epsilon_1$  such that
$\epsilon_1<2^{-2}$.

\medskip

Assume that we have constructed/chosen the following objects with
the listed properties:

\medskip

\begin{itemize}

\item[(1)] Smoothly bounded simply connected domains $U_j\subset\triangle$ and
$\mathfrak{b}$-nice domains $W_j$ such that $U_1\subset\subset
U_2\subset\subset\cdot\cdot\cdot\subset\subset U_N\subset\subset
W_N\subset W_{N-1}\subset\cdot\cdot\cdot\subset
W_1\subset\overline\triangle$,

\item[(2)] Embeddings $\phi_j:W_j\emb\mathbb C^2$, such that
$\phi_j(U_j)\subset\mathbb B_j$ and $\|\phi_j(z)\|=j$ for all
$z\in\partial U_j$,

\item[(3)] $\phi_j(U_j\setminus U_{j-1})\subset\mathbb
C^2\setminus\overline{\mathbb B}_{j-1-2^{-j}}$,

\item[(4)] $\phi_j(\overline U_j)\cap X=\emptyset$.

\item[(5)] The pair $(\phi_N,W_N)$ satisfies the condition in Lemma
\ref{lemma}.

\end{itemize}

\

Moreover we have inductively chosen a sequence
$\epsilon_1>\epsilon_2>\cdot\cdot\cdot >\epsilon_N$  with
$\epsilon_j<2^{-j-1}$ and assured that

\

\begin{itemize}

\item[(6)] If $\psi:\overline U_j\rightarrow\mathbb C^2$ is a
holomorphic map with $\|\psi-\phi_j\|_{\overline{U_j}}<\epsilon_j$
then $\psi$ is an embedding and $\psi(\overline U_j)\cap
X=\emptyset$,

\item[(7)]
$\|\phi_j-\phi_{j-1}\|_{\overline{U_{j-1}}}<\epsilon_{j-1}2^{-j}$.

\end{itemize}

\

We now show how to get $U_{N+1},\phi_{N+1},W_{N+1}$ and
$\epsilon_{N+1}$  so that we have $(1)-(7)$ with $N+1$ in place of
$N$.

\

To apply Lemma \ref{lemma} we let $\phi:=\phi_N, W:=W_N,
V:=\overline{U}_N$ and $\epsilon:=\epsilon_N 2^{-N-1}$. (Technically
$W_N$ is not the same as in the lemma, but by the Riemann mapping
theorem it does not make a difference).  Let $\phi_{N+1}$ and
$W_{N+1}$ denote the objects corresponding to $\tilde\phi$ and
$\tilde W$ in the conclusion of the lemma.  We get immediately then
that the pair $(\phi_{N+1},W_{N+1})$ satisfy the conditions in the
lemma, i.e. we have $(5)$.  In particular this means that
$$
\phi_{N+1}(W_{N+1})\cap(\overline{\mathbb B}_{N+1}\cap X)=\emptyset.
$$
Since $\overline U_N = V$ by assumption we also get that
$\|\phi_{N+1}-\phi_N\|_{\overline U_N}<\epsilon=\epsilon_N2^{-N-1}$,
i.e. we get $(7)$.  To define $U_{N+1}$ we consider the set
$S:=\phi_{N+1}^{-1}(\phi_{N+1}(W_{N+1})\cap\overline{\mathbb
B}_{N+1})$. Note first that $U_N\subset S$ since we just established
$(7)$.  This means that we may define $U_{N+1}$ to be the interior of the connected
component of $S$ that contains $U_N$.  By $(e)$ we have that
$U_{N+1}$ is smoothly bounded, and we get $(1),(2)$ and $(4)$. Since
$U_{N+1}\setminus U_N\subset W_{N+1}\setminus U_N$ we get $(3)$ from
$(g)$.  Finally we choose $\epsilon_{N+1}$ small enough to get
$(6)$.

\medskip

To finish the proof we construct a sequence $(U_j,\phi_j)$ according
to the above procedure.  We define $U:=\cup_{j=1}^\infty U_j$.  Then
$U$ is an increasing union of simply connected domains and so $U$ is
itself simply connected.  By the Riemann Mapping Theorem $U$ is
conformally equivalent to the unit disk.  We define a map
$\psi:U\rightarrow\mathbb C^2$ by
$$
\psi(z)=\lim_{j\rightarrow\infty}\phi_j(z).
$$
To see that this is well defined we consider a point $z\in\overline
U_k$: For $m>n\geq k$ we have by $(7)$ that
$$
\|\phi_m(z)- \phi_n(z)\|\leq\sum_{i=n+1}^m
\|\phi_i(z)-\phi_{i-1}(z)\|\leq\sum_{i=n+1}^m\epsilon_{i-1}2^{-i}<\epsilon_n<2^{-n-1}.
$$
This shows that $\{\phi_j(z)\}$ is a Cauchy-sequence and so $\psi$
defines a holomorphic map from $U$ into $\mathbb C^2$.  It also
shows that
$$
(*) \ \|\psi(z)-\phi_k(z)\|\leq\epsilon_k.
$$
for all $z\in\overline U_k$, so it follows from $(6)$  that $\psi$
is an embedding, and that $\psi(U)\cap X=\emptyset$.   \

To see that $\psi$ is proper, consider a point $z\in
U_{k+1}\setminus\overline U_k$ for some $k$.  By (3) we have that
$\|\phi_{k+1}(z)\|\geq k-2^{-k-1}$, and so by $(*)$ we get that
$\|\psi(z)\|\geq k-2^{-k-1}-\epsilon_k>k-2^{-k}$. This means that
$\|\psi(z)\|>k-2^{-k}$ for all $z\in U\setminus U_k$, hence $\psi$
is proper. $\hfill\square$

\section{Concluding remarks}

We note that if $n> 2$ the analogous result to our main theorem holds with $\C^n$
instead of $\C^2$ ( if the codimension of $X$ is bigger than one, this is an easy consequence
of transversality).

\begin{theorem}
Let $X$ be any closed proper subvariety of $\mathbb C^n, n \ge 2$ and let
$\triangle$ denote the unit disk in $\mathbb C$.  Then there exists
a proper holomorphic embedding $\phi:\triangle\emb\mathbb C^n$ such
that $\phi(\triangle)\cap X=\emptyset$.
\end{theorem}

The proof of Theorem \ref{main} works in this case also, but one could simply take a 
$2$-dimensional complex-linear subspace $C$ of $\mathbb C^n$ such that
$X\cap C\neq C$, and embed $\triangle$ into $C\setminus(X\cap C)$ according to the 
main theorem. \

As pointed out in the introduction there is no way to add interpolation conditions, but we believe 
that it is possible to make the image of the embedding containing a prescribed discrete  subset $A$ of $\C^n$ (with $A\cap X = \emptyset$).

\bibliographystyle{amsplain}

 \end{document}